\documentclass[a4,11pt]{article}

\usepackage{amssymb}
\usepackage{graphicx}

\newtheorem{thm}{Theorem}[section]
\newtheorem{rem}[thm]{Remark}
\newtheorem{lem}[thm]{Lemma}

\title{Vertex unfoldings of tight polyhedra}

\author{Toshiki Endo\thanks{Jiyu Gakuen College, {\tt end@prf.jiyu.ac.jp}} \and Yuki Suzuki\thanks{Faculty of Science Division II, Tokyo University of Science, {\tt j2111076@ed.tus.ac.jp}}}
\date{}

\begin{document}

\maketitle

\begin{abstract}
The unfolding of a polyhedron along its edges is known as a vertex unfolding if adjacent faces are allowed to be connected not only at an edge but also at a vertex. 
Demaine et al. \cite{DEHO} showed that every triangulated polyhedron has a vertex unfolding. 
We extend this result to a tight polyhedron, where a polyhedron is tight if its non-triangular faces are mutually non-incident. 
\end{abstract}

%% \begin{keyword}
%% polyhedron \sep edge unfolding \sep vertex unfolding \sep planar graph \sep vertex-face tour
%% \end{keyword}

\section{Introduction}
We investigate a procedure to cut open a polyhedron along its edges and unfold it to a connected flat piece without overlap. 
The unfolding needs to consist of the faces of the polyhedron joined along the edges. 
This type of unfolding has been referred to as {\it edge unfolding} or simply {\it unfolding}. 
It is known that some non-convex polyhedra have no edge unfoldings. 
However, no example of a convex polyhedron that has no edge unfolding is known. 
The determination of whether every convex polyhedron has an edge unfolding is a long-standing open problem. 
The difficulty of this question led to the exploration of other unfoldings that have a broader definition of edge unfolding. 
We pay attention to a {\it vertex unfolding} that permits two faces joined not only at an edge but also at a vertex, that is, the resulting piece may have a disconnected interior. 
See \cite[\S 22]{DO} for details of edge unfolding and vertex unfolding. 

In \cite{DEHO}, Demaine et al. showed the following, where they proved conclusively that $\cal P$ does not need to be a polyhedron, but may be a connected triangulated 2-manifold, possibly with boundaries. 

\begin{thm}
\emph{(Demaine et al \cite{DEHO})}
\label{thm:Demaine}
Let $\cal P$ be a polyhedron. 
If $\cal P$ is triangulated, then it has a vertex unfolding. 
\end{thm}

We broadly describe the proof of Theorem \ref{thm:Demaine} here and describe it in detail in the following sections. 
Their algorithm \cite{DEHO} first finds a spanning path from triangle to triangle on the surface of the polyhedron, connecting through common vertices, and further, lays out the triangles along a line without overlap. 

Their method is based on the condition that all faces are triangular, and the existences of the face path and the line-layout of it might actually fail for a polyhedron with non-triangular faces. 
For example, the truncated cube has no face path since 
its six octagons are inadequate to lay out eight triangles along a line, and if a face path consists of isosceles trapezoids, a local overlap might occur in a long strip. 

In this paper, we fix these problems and make progress on Theorem \ref{thm:Demaine} for a polyhedron with non-triangular faces. 
A (possibly non-convex) polyhedron $\cal P$ is {\it tight} if no two non-triangular faces share a vertex. 
Examples of tight polyhedra are the snub cube, snub dodecahedron, pyramids, and antiprisms. 
The main theorem in this paper is as follows.

\begin{thm}
\label{thm:main}
Let $\cal P$ be a polyhedron. 
If $\cal P$ is tight, then it has a vertex unfolding. 
\end{thm}

\begin{figure}[htb]
\centerline{
\includegraphics[width=0.8\textwidth]{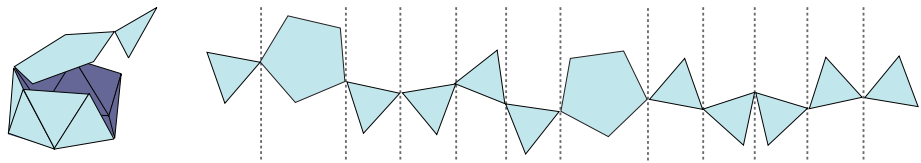}}
\caption{Vertex unfolding of the pentagonal antiprism}
\label{fig:layout}
\end{figure}

Figure \ref{fig:layout} shows a vertex unfolding of the pentagonal antiprism. 
Our proof basically depend on the method in \cite{DEHO}. 
We will describe it in Sections \ref{sec:Hamiltonian vertex-face tour}-\ref{sec:arrangement of the faces}. 
Our new result is a graph theoretical part of it, which is contained in Section \ref{sec:Hamiltonian vertex-face tour}.  

\section{Hamiltonian vertex-face tour}
\label{sec:Hamiltonian vertex-face tour}
In this section, we observe tight polyhedra from a graph theoretical standpoint. 
We use standard terminology and notations of graph theory, for examples, see \cite{D}. 
By Steinitz's theorem, a surface of a polyhedron corresponds to a 3-connected plane graph. 
Thus, we also call a 3-connected plane graph {\it tight} if its non-triangular faces are mutually non-incident. 
We prepare some more definitions. 

Let $G$ be a tight graph. 
A disjoint union $T$ of closed alternating sequences of vertices $v_i$ and faces $f_i$ of $G$ is called a {\it spanning vertex-face tour} if each face of $G$ appears exactly once in $T$ and each closed component $(v_1$, $f_1$, $v_2$, $f_2$, $\cdots$, $v_{k}$, $f_{k}$, $v_1)$ satisfies that $v_i$ and $v_{i+1}$ are distinct and both are incident to the face $f_{i}$ for $i=1, 2, \cdots, k$ (indices are taken modulo $k$). 
Some vertex of $G$ may be repeated in $T$; conversely, some vertex may not appear in $T$. 
If a spanning vertex-face tour $T$ is connected, then $T$ is called a {\it Hamiltonian vertex-face tour}. 

\begin{figure}[htp]
\centerline{
\includegraphics[width=0.8\textwidth]{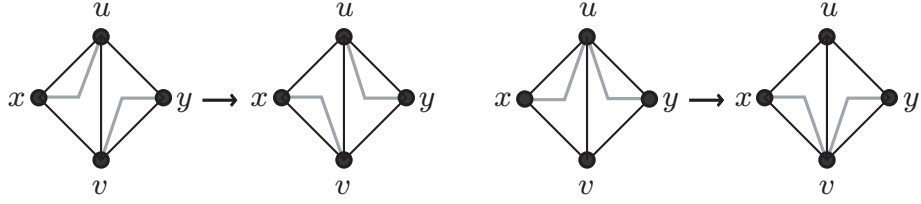}}
\caption{Switching operation (left) and reflecting operation (right)}
\label{fig:operation}
\end{figure}

Next, we define two operations on a spanning vertex-face tour $T$. 
Let $f=uvx$ and $f'=uvy$ be two adjacent triangular faces of $G$. 
We refer to the operation of replacing $(u, f, x)$ and $(v, f', y)$ with $(v, f, x)$ and $(u, f', y)$, respectively, as the {\it switching} operation, and the operation of replacing $(u, f, x)$ and $(u, f', y)$ with $(v, f, x)$ and $(v, f', y)$, respectively, as the {\it reflecting} operation (Figure \ref{fig:operation}). 
Note that simultaneous changing of combinations between vertices and faces at several triangular faces of $T$, such as in a switching operation or reflecting operation, may produce another spanning vertex-face tour $T'$. 
In general, we refer to such operations as {\it triangular recombinations}. 

First, we prove the following lemma.

\begin{lem}
\label{lem:changing diagonal}
Let $G$ be a tight graph. 
Let $F^{\ast}=\{f_{1}, f_{2}, \cdots, f_{m}\}$ be the set of all non-triangular faces of $G$, and let $v_{i}$ and $v'_{i}$ be distinct vertices of $f_{i}$ for $i=1,2,\cdots, m$. If $G$ has a spanning vertex-face tour $T$, then $G$ has a spanning vertex-face tour $T'$ containing each $(v_{i}$, $f_{i}$, $v'_{i})$ for $i=1,2,\cdots, m$. 
\end{lem}

\noindent
{\sc Proof.}\quad
For simplicity, let $f=v_1v_2\cdots v_n$ denotes any non-triangular face of $G$.  Let $g_i$ be the triangular face adjacent to $f$ by sharing $v_iv_{i+1}$ for $i=1,2,\cdots n$ (indices are taken modulo $n$), and let $u_{i}$ be the remaining vertex of $g_i$ for $i=1,2,\cdots n$. 

We only have to show that if $T$ contains $(v_1, f, v_{k})$ for some $2 \leq k\leq n-1$, then $T$ can be converted to a spanning vertex-face tour $T'$ containing $(v_1$, $f$, $v_{k+1})$ instead of $(v_1, f, v_{k})$ by performing only triangular recombinations. 
\medskip

\noindent
Case 1: $T$ contains $(u_k, g_k, v_{k})$ or $(u_k, g_k, v_{k+1})$. 

In this case, we can obtain $T'$ from $T$ by replacing $(v_1, f, v_{k})$ with $(v_1$, $f$, $v_{k+1})$, and simultaneously by replacing $(u_k$, $g_k$, $v_{k})$ with $(u_k$, $g_k$, $v_{k+1})$ in the former case and $(u_k$, $g_k$, $v_{k+1})$ with $(u_k$, $g_k$, $v_{k})$ in the latter case. 
\medskip

\noindent
Case 2: $T$ contains $(v_k, g_k, v_{k+1})$. 

In this case, we check the triangular faces incident to $v_i$ from $i=k+1, k+2,\cdots, n-1, n, 1,2,\cdots, k$ in turn. 
For $i=k+1, k+2,\cdots, n-1, n, 1,2,\cdots, k$, if they exist, let $h_{i}^1, \cdots, h_{i}^{p_i-3}$ be the triangular faces incident to $v_{i}$ between $g_{i-1}$ and $g_{i}$ and opposite to $f$ in cyclic order, where $p_i=\deg v_{i}$, and let $w_{i}^{1}, w_{i}^{2}, \cdots, w_{i}^{p_i-4}$ be the vertices incident to $v_i$ from $u_{i-1}$ to $u_i$. 

First, we examine the triangular faces incident to $v_{k+1}$. 
If $T$ contains $(u_{k}$, $h_{k+1}^1$, $w_{k+1}^1)$, then the switching operation at $g_{k}$ and $h_{k+1}^1$ leads this case to Case 1. 
If $T$ contains $(v_{k+1}, h_{k+1}^1, w_{k+1}^1)$, then the reflecting operation at $g_{k}$ and $h_{k+1}^1$ again leads this case to Case 1. 
Thus, $T$ must contain $(v_{k+1}, h_{k+1}^1, u_{k})$. 
By repeating this argument, we can say that $T$ contains $(v_{k+1}, g_{k+1}, u_{k+1})$. 

Second, we check the triangular faces incident to $v_{k+2}$, and we can say that $T$ contains $(v_{k+2}, g_{k+2}, u_{k+2})$. 
Repeating this argument, we finally deduce that $T$ contains $(v_{k}, h_{k}^{p_{k}-3}, w_{k}^{p_{k}-4})$. 
Thus, we can apply the reflecting operation at $h_{k}^{p_{k}-3}$ and $g_k$, which leads this case to Case 1. 
\hfill$\Box$
\bigskip

\begin{rem}
\normalfont
\label{rem:triangle}
In the proof of Lemma \ref{lem:changing diagonal}, we can choose a triangular face $f$ as a member of $F^{\ast}$ if the faces incident to $f$ are all triangular faces. 
\end{rem}

Next, we prove the following lemma.

\begin{lem}
\label{lem:triangulation}
Let $G$ be a plane triangulation. 
Then $G$ has a spanning vertex-face tour $T$. 
\end{lem}

\begin{figure}[htb]
\centerline{
\includegraphics[width=0.5\textwidth]{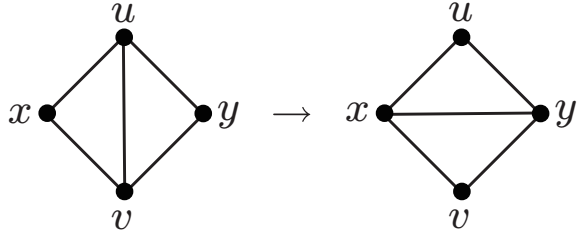}}
\caption{Diagonal flip}
\label{fig:flip}
\end{figure}

\begin{figure}[htb]
\centerline{
\includegraphics[width=0.45\textwidth]{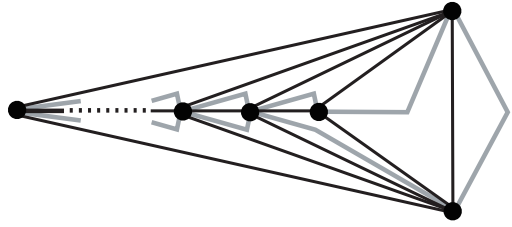}}
\caption{Standard triangulation and Hamiltonian vertex-face tour (gray line)}
\label{fig:standard}
\end{figure}

In order to prove Lemma \ref{lem:triangulation}, we use Wagner's theorem \cite{W}, which states that every triangulation can be transformed into the standard triangulation by a finite sequence of diagonal flips. 
Here, the operation {\it diagonal flip} is defined as follows. 
Let $uv$ be an edge of a triangulation $G$. 
Let $uvx$ and $uvy$ be the faces incident to $uv$. 
Then $x$ and $y$ are distinct vertices unless $=K_3$. 
If $x$ and $y$ are not adjacent, then a diagonal flip is performed to obtain a new triangulation $G'$ from $G$ by deleting $uv$ and adding the edge $xy$ (Figure \ref{fig:flip}). 
The standard triangulation is defined as illustrated in Figure \ref{fig:standard}. Note that the standard triangulation has a Hamiltonian vertex-face tour. 

\noindent
{\sc Proof.}\quad
Let $f_1=v_1v_2v_3$ and $f_2=v_3v_4v_1$ be two adjacent faces of $G$, $G'$ be the triangulation obtained from $G$ by performing the diagonal flip at $v_1v_3$, and $f_1'=v_2v_3v_4$ and $f_2'=v_4v_1v_2$ be the new faces of $G'$. 
From Wagner's theorem and the fact in Figure \ref{fig:standard}, we only have to show that if $G$ has a spanning vertex-face tour $T$ then $G'$ has a spanning vertex-face tour $T'$. 

Let $g_1$ and $g_2$ be triangular faces of $G$ that are adjacent to $f_1$ by sharing $v_1v_2$ and $v_2v_3$, respectively, let $g_3$ and $g_4$ be triangular faces of $G$ that are adjacent to $f_2$ by sharing $v_3v_4$ and $v_4v_1$, respectively, and let $u_{i}$ be the remaining vertex of $g_i$ for $i=1, 2, 3, 4$. 

If they exist, let $h_i^1, h_i^2,\cdots, h_i^{p_i-3}$ for $i=2, 4$ and $h_i^1, h_i^2,\cdots, h_i^{p_i-4}$ for $i=1, 3$ be the triangular faces that are incident to $v_i$ and between $g_{i-1}$ and $g_{i}$ but opposite to $f_1$ and $f_2$ in cyclic order, where $p_i=\deg_Gv_{i}$. 
Let $w_{i}^{1}, w_{i}^{2}, \cdots, w_{i}^{p_i-4}$ for $i=2, 4$ and $w_{i}^{1}, w_{i}^{2}, \cdots, w_{i}^{p_i-5}$ for $i=1, 3$ be the vertices incident to $v_i$ from $u_{i-1}$ to $u_i$. 
We divide the proof into four cases. 
\medskip

\noindent
Case 1: $T$ contains $(v_2, f_1, v_3)$ and $(v_4, f_2, v_1)$. 

In this case, we can obtain $T'$ from $T$ by replacing $(v_2, f_1, v_3)$ and $(v_4, f_2, v_1)$ with $(v_2, f_1', v_3)$ and $(v_4, f_2', v_1)$, respectively. 
\medskip

\noindent
Case 2: $T$ contains $(v_1, f_1, v_3)$ and $(v_4, f_2, v_1)$. 

We consider the triangular faces incident to $v_2$ from $g_1$ to $g_2$ in turn. 
If $T$ contains $(u_1, g_1, v_2)$, then the switching operation at $f_1$ and $g_1$ leads the case to Case 1. 
If $T$ contains $(u_1, g_1, v_1)$, then the reflecting operation at $f_1$ and $g_1$ again leads the case to Case 1. 
Thus, $T$ must contain $(v_2, g_1, v_1)$. 
Next, we consider $h_2^{1}$ and similarly we can say that $T$ must contain $(v_2, h_2^{1}, u_1)$; by repeating this argument until we reach $g_2$, we can say that $T$ contains $(v_2, g_2, u_2)$. 

Further, we examine the triangular faces incident to $v_3$, and similarly we can say that $T$ contains $(v_3, h_3^{1}, u_2), (v_3, h_3^{2}, w_3^1), \cdots, (v_3, g_3, u_3)$. 
Repeating this argument for $v_3, v_4$, and finally $v_1$ in turn, we can say that $T$ contains $(v_1, h_1^{p_1-4}, w_1^{p_1-5})$.
Thus, we can apply the reflecting operation at $g_1$ and $h_1^{p_1-4}$. 
\medskip

\noindent
Case 3: $T$ contains $(v_1, f_1, v_2)$ and $(v_4, f_2, v_1)$. 

In this case, we consider the triangular faces incident to $v_3$ from $g_2$ to $g_3$ in turn. 
If $T$ contains $(u_2, g_2, v_3)$, then the switching operation at $f_1$ and $g_2$ leads the case to Case 2. 
If $T$ contains $(u_2, g_2, v_2)$, then the reflecting operation at $f_1$ and $g_2$ again leads the case to Case 2. 
Thus, $T$ must contain $(v_3, g_2, v_2)$. 
Similarly, we can say that $T$ contains $(v_3, h_3^{1}, u_2), (v_3, h_3^{2}, w_3^1), \cdots, (v_3, g_3, u_3)$. 
Thus, if we perform the switching operation at $g_3$ and $f_2$, the situation becomes a symmetric version of Case 2. 
\medskip

\noindent
Case 4: $T$ contains $(v_1, f_1, v_3)$ and $(v_3, f_2, v_1)$. 

In this case, we consider the triangular faces incident to $v_3$ from $g_2$ to $g_3$ in turn. 
If $T$ contains $(u_2, g_2, v_2)$, then the switching operation at $f_1$ and $g_2$ leads the case to the symmetric version of Case 2. 
If $T$ contains $(u_2, g_2, v_3)$, then the reflecting operation at $f_1$ and $g_2$ again leads the case to the symmetric version of Case 2. 
Thus, $T$ must contain $(v_3, g_2, v_2)$. 
Similarly, we can say that $T$ contains $(v_3, h_3^{1}, u_2), (v_3, h_3^{2}, w_3^1), \cdots, (v_3, g_3, u_3)$. 
Thus, if we apply the reflecting operation at $g_3$ and $f_2$, the situation becomes the same as Case 2. 
\hfill$\Box$
\bigskip

\begin{lem}
\label{lem:spanning vertex-face tour}
Let $G$ be a tight graph. 
Then $G$ has a spanning vertex-face tour $T$. 
\end{lem}

\noindent
{\sc Proof.}\quad
We prove this by applying a double-induction on the size and the number of the maximum face of $G$. 
\medskip

\noindent
Case 1: $G$ has no non-triangular faces. 

This case follows from Lemma \ref{lem:triangulation}. 
\medskip

\noindent
Case 2: The maximum face size of $G$ is at least four. 

Let $f=v_1v_2 \cdots v_n$ be a face with the maximum size ($n\geq 4$). 
From the planarity of $G$, we may assume that $G'=G+v_1v_3$ is tight. 
Let $f'=v_1v_2v_3$ and $f''=v_3v_4\cdots v_nv_1$ be the new faces of $G'$. 
From the inductive hypothesis, $G'$ has a spanning vertex-face tour $T'$. 
We show that $G$ has a spanning vertex-face tour $T$. 

From Lemma \ref{lem:changing diagonal} and Remark \ref{rem:triangle}, we may assume that $T'$ contains $(v_1, f'', v_3)$. 
Therefore, if $T'$ contains $(v_2, f', v_1)$, then we can obtain $T$ from $T'$ by replacing $(v_2, f', v_1)$ and $(v_1, f'', v_3)$ with $(v_2, f, v_3)$; the case where $T'$ contains $(v_2, f', v_3)$ is similar. 
Thus, we may assume that $T'$ contains $(v_1, f', v_3)$. 

We consider the triangular faces incident to $v_i$ from $i=3, 4, \cdots, n$ in turn. 
Let $g_2$ be the triangular face adjacent to $f'$ sharing $v_2v_{3}$, and 
for $i=3, 4, \cdots, n$, let $g_i$ be the triangular face adjacent to $f''$ sharing $v_iv_{i+1}$ (indices are taken modulo $n$). 
Further, let $u_{i}$ be the remaining vertex of $g_i$ for $i=2,3,\cdots n$. 

First, we examine the triangular faces incident to $v_3$ from $g_2$ to $g_3$. 
If they exist, let $h_3^1$, $h_3^2$, $\cdots$, $h_3^{p_3-4}$ be the triangular faces between $g_2$ and $g_3$ in cyclic order, where $p_3=\deg_{G'}v_3$, and let $h_{3}^1=v_{3}u_{2}w_{3}^{1}$, $h_{3}^j=v_{3}w_{3}^{j-1}w_{3}^{j}$ for $j=2,3,\cdots, p-6$, and $h_{3}^{p-4}=v_{3}w_{3}^{p-5}u_3$. 
If $T'$ contains $(v_2, g_{2}, u_2)$, then the switching operation at $f'$ and $g_{2}$ yields a new spanning vertex-face tour containing $(v_1, f', v_2)$, and if $T'$ contains $(v_{3}, g_{2}, u_2)$, then the reflecting operation at $f'$ and $g_{2}$ yields a new vertex-face tour containing $(v_1, f', v_2)$, and in both cases, we can obtain $T$ by replacing $(v_1, f', v_2)$ and $(v_1, f'', v_3)$ with $(v_2, f, v_3)$. 
Thus, $T'$ must contain $(v_{3}, g_{2}, v_2)$. 
By repeating this argument from $ h_{3}^1$ to $g_{3}$, we can say that $T'$ contains $(v_{3}, g_{3}, u_{3})$. 
Thus, we can obtain a new vertex-face tour by replacing $(v_1, f'', v_3)$ and $(v_{3}, g_{3}, u_{3})$ with $(v_{1}, f'', v_4)$ and $(v_4, g_3, u_3)$, respectively. 
In this case, we can obtain $T$ by replacing $(v_1, f', v_3)$, $(v_1, f'', v_4)$ with $(v_3, f, v_4)$. 
\hfill$\Box$
\bigskip

\begin{lem}
\label{lem:Hamiltonian vertex-face tour}
Let $G$ be a tight graph. 
Let $F^{\ast}=\{f_{1}, f_{2}, \cdots, f_{m}\}$ be the set of all non-triangular faces of $G$, and let $v_{i}$ and $v'_{i}$ be distinct vertices of $f_{i}$ for $i=1,2,\cdots, m$. 
If $G$ has a spanning vertex-face tour $T'$ containing each $(v_{i}$, $f_{i}$, $v'_{i})$ for $i=1,2,\cdots, m$, then $G$ has a Hamiltonian vertex-face tour $T''$ containing each $(v_{i}$, $f_{i}$, $v'_{i})$ for $i=1,2,\cdots, m$. 
\end{lem}

\noindent
{\sc Proof.}\quad
We show that $T'$ can be converted to be a connected spanning vertex-face tour by performing only a series of triangular recombinations. 
Suppose that $T'$ is disconnected at two adjacent faces $g_1$ and $f_1$. 
We may assume that $g_1=v_1v_2u_1$ is a triangular face. 
Let $f_1=v_1v_2\cdots v_n$. 
We divide the proof into two cases. 
\medskip

\noindent
Case 1: $n=3$. 

In this case, we may assume that two components of $T'$ containing $(v_1, g_1, u_1)$ and $(v_2, f_1, v_3)$ are disconnected. 
Thus, we can make the two components connected by performing the switching operation at $g_1$ and $f_1$. 
\medskip

\noindent
Case 2: $n\geq 4$. 

Suppose that $T'$ contains $(v_{k_1}, f_1, v_{k_2})$ for some $k_1$ and $k_2$. 
Then $T'$ must contain $(v_{k_1}, h_{k_1}^l, w_l)$ for some triangular face $h_{k_1}^l$ incident to $v_{k_1}$, and for some vertex $w_l$ of $h_{k_1}^l$. 
Now, $h_{k_1}^l$ is connected to $g_1$ by a path of triangular faces, and it holds from Case 1 that they can become connected by triangular recombinations. 
\hfill$\Box$
\bigskip

Our goal is the following. 

\begin{thm}
\label{thm:Hamiltonian vertex-face tour}
Let $G$ be a tight graph. 
Let $F^{\ast}=\{f_{1}, f_{2}, \cdots, f_{m}\}$ be the set of all non-triangular faces of $G$, and let $v_{i}$ and $v'_{i}$ be distinct vertices of $f_{i}$ for $i=1,2,\cdots, m$. 
Then $G$ has a Hamiltonian vertex-face tour containing each $(v_{i}$, $f_{i}$, $v'_{i})$ for $i=1,2,\cdots, m$. 
\end{thm}

\noindent
{\sc Proof.}\quad
Let $G$ be a tight graph. 
From Lemma \ref{lem:spanning vertex-face tour}, $G$ has a spanning vertex-face tour $T$. 
Then, from Lemma \ref{lem:changing diagonal}, $G$ has a spanning vertex-face tour $T'$ containing each $(v_{i}$, $f_{i}$, $v'_{i})$. 
Thus, from Lemma \ref{lem:Hamiltonian vertex-face tour}, $G$ has a Hamiltonian vertex-face tour $T''$ containing each $(v_{i}$, $f_{i}$, $v'_{i})$. 
\hfill$\Box$
\bigskip

\section{Non-crossing Hamiltonian face path}
\label{sec:non-crossin}
For a polyhedron $\cal P$ and its graph $G$, a Hamiltonian vertex-face tour of $G$ guarantees an existence of a path of the faces of $\cal P$. 
We call it a {\it Hamiltonian face path} of $\cal P$. 
However, the path might cross itself in the sense that it contains the pattern $(\cdots, f_1, v, f_3, \cdots, f_2, v, f_4, \cdots)$ with the faces $f_1, f_2, f_3, f_4$ incident to a vertex $v$ appearing in cyclic order. 
This make it physically impossible for the faces of an unfolding to be a single piece. 
Hence, we need to detect a non-crossing path. 
A face path of $\cal P$ (likewise, a vertex-face tour of $G$) is {\it non-crossing} if it has no patterns as that described above.  

\begin{lem}
\label{lem:non-crossing Hamiltonian vertex-face tour}
In Theorem $\ref{thm:Hamiltonian vertex-face tour}$, any Hamiltonian vertex-face tour of $G$ can be converted to a non-crossing one. 
\end{lem}

\begin{figure}[htb]
\centerline{
\includegraphics[width=0.6\textwidth]{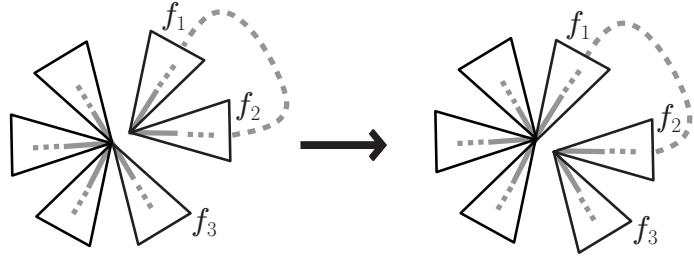}}
\caption{Converting a Hamiltonian vertex-face tour into a non-crossing one}
\label{fig:noncrossing}
\end{figure}

\noindent
{\sc Proof.}\quad
This is contained in \cite{DEHO}. 
The key point of the proof is as follows. 
Suppose that a Hamiltonian vertex-face tour $T$ crosses at a vertex $v$. 
Let $f_1$, $f_2$, $\cdots$ be the faces passing through $v$ in $T$ in cyclic order. 
We remove the face path of $(\cdots, f_1, v, f_2, \cdots)$ from $T$ as depicted in Figure \ref{fig:noncrossing}. 
If the resulting tour is disconnected, then we remove the face path of $(\cdots, f_2, v, f_3, \cdots)$ instead of
the above from $T$ such that the resulting tour is connected.  
By repeating this operation at every vertex of $G$, we obtain a non-crossing Hamiltonian vertex-face tour. 
\hfill$\Box$
\bigskip

\section{Layout of a face path}
\label{sec:arrangement of the faces}
In this section, we exhibit the procedure to lay out the faces of a tight polyhedron $\cal P$ to form a vertex unfolding. 
First, we show the following. 

\begin{lem}
\label{lem:polygon}
Let $\cal F$ be a (possibly non-convex) polygon with four sides or more. 
Then, there are two vertices $u, v$ of $\cal F$ and an arrangement of $\cal F$ in a vertical interval of the plane with $u$ and $v$ on its left and right boundaries, respectively. 
\end{lem}

\noindent
{\sc Proof.}\quad
We only have to choose $u$ and $v$ such that the length of segment $\overline{uv}$ is longest among all diagonals and edges of $\cal F$. 
\hfill$\Box$
\bigskip

\noindent
{\sc Proof of Theorem\ \ref{thm:main}.}\quad
Let $\cal P$ be a tight polyhedron. 
Consider the graph $G$ of $\cal P$. 
From Lemma \ref{lem:non-crossing Hamiltonian vertex-face tour}, $G$ has a non-crossing Hamiltonian vertex-face tour $T$. 
Let $\cal T$ be the corresponding face path of $\cal P$. 
We may assume from Theorem \ref{thm:Hamiltonian vertex-face tour} that $\cal T$ uses the vertices of Lemma \ref{lem:polygon} in each non-triangular face. 

Now, we can arrange the faces as follows; this is a consequence of Lemma 22.6.2 in textbook \cite{DO}. 
Suppose inductively that $\cal P$ has been laid out along a line up to face $f_{i-1}$ with all faces left of vertex $v_i$, which is the rightmost vertex of $f_{i-1}$. 
Let $(v_i, f_i, v_{i+1})$ be the next face in $\cal T$. 
If $f_i$ is a triangular face, rotate $f_i$ around $v_i$ such that $f_i$ lies horizontally between or at the same horizontal coordinate as $v_i$ and $v_{i+1}$. 
If $f_i$ is a non-triangular face, we can use Lemma \ref{lem:polygon}. 
Repeating this process along $\cal T$ produces a non-overlapping layout of the faces of $\cal P$. 
Thus, $\cal P$ has a vertex unfolding. 
\hfill$\Box$


\bigskip

\begin{thebibliography}{0}

\bibitem[1]{DEHO}
E. D. Demaine, D. Eppstein, G. W. Hart, and J. O'Rourke, Vertex-unfoldings of simplicial manifolds, in Discrete Geometry (ed. Andras Bezdek), Marcel Dekker, New York, 2003, 215--228. 

\bibitem[2]{DO}
E. D. Demaine and J. O'Rourke, Geometric Folding Algorithms: Linkages, Origami, Polyhedra,  Cambridge University Press, Cambridge (2007). 

\bibitem[3]{D}
R. Diestel, Graph Theory (3e), Graduate Texts in Mathematics No. 173, Springer (2005).

\bibitem[4]{W}
K. Wagner, Bemerkung zum Vierfarbenproblem, Jber. Deutsch. Math.-Verein. 46 (1936), 26--32.
\end{thebibliography}
\end{document}